\documentclass[12pt]{amsart}
\usepackage{amsmath, amssymb}
\usepackage{hyperref}

\author{}
\title{}
\date{}

\newcommand{\bbN}{{\mathbb N}}
\newcommand{\bbQ}{{\mathbb Q}}
\newcommand{\bbR}{{\mathbb R}}
\newcommand{\bbZ}{{\mathbb Z}}

\newtheorem{lemma}{Lemma}[section]

\newtheorem{theorem}{Theorem}

\begin{document}
\baselineskip=17pt

\begin{center}
{\bf COMPLETE SOLUTION OF THE DIOPHANTINE EQUATION} $X^{2}+1=dY^{4}$ 
{\bf AND A RELATED FAMILY OF QUARTIC THUE EQUATIONS} \\

\vspace{5.0mm} 

{\sc Chen Jian Hua} 

\vspace{4.0mm} 

{\em The Electric Power Testing and Research Institute of Hubei Province} \\ 
{\em Wuhan 430077, P. R. China} 

\vspace{5.0mm} 

AND 

\vspace{5.0mm} 

{\sc Paul M. Voutier} 

\vspace{4.0mm} 

{\em Department of Mathematics, City University} \\
{\em Northampton Square, London EC1V 0HB, UK}
\end{center}

\vspace{7.0mm}

\begin{quote}
{\footnotesize In this paper, we use the method of Thue and Siegel, 
based on explicit Pade approximations to algebraic functions, to 
completely solve a family of quartic Thue equations. From this result, 
we can also solve the diophantine equation in the title. We prove 
that this equation has at most one solution in positive integers 
when $d \geq 3$. Moreover, when such a solution exists, it is of 
the form $(u,\sqrt{v})$ where $(u,v)$ is the fundamental solution 
of $X^{2}+1=dY^{2}$.} 
\end{quote}

\vspace{5.0mm}

\begin{center}
1. {\sc Introduction}
\end{center}

\setcounter{section}{1}

\vspace{3.0mm}

The goal of this paper is the solution of the diophantine equation 
\begin{equation}
\label{eq:1.1}
X^{2} + 1 = dY^{4}.  
\end{equation}

Many people have studied this equation throughout the history of 
diophantine equations. Previously, the best known result was 

\begin{theorem}[Ljunggren \cite{Ljung}]
\label{thm:A}
If the fundamental unit of the quadratic field $\bbQ (\sqrt{d})$ is 
not also the fundamental unit of the ring $\bbZ [\sqrt{d}]$, then  
equation $(\ref{eq:1.1})$ has at most two possible solutions in 
positive integers, and these can be found by a finite algorithm.
\end{theorem}

His proof is exceedingly complicated (see Theorem 10 and the 
surrounding text on p.271 of \cite{Mord}). 

In another article \cite{Chen2}, the first author proved, that when 
$d>0$ is large enough, there exists only one possible solution. More 
precisely, we have 

\begin{theorem}[Chen \cite{Chen2}] 
\label{thm:B}
Let $d>0$ and put $\epsilon = u+v\sqrt{d}$ where $(u,v)$ 
is the least positive solution of the Pell equation 
\begin{displaymath}
X^{2} + 1 = dY^{2}.  
\end{displaymath}

If $\epsilon > 5 \times 10^{7}$, then the equation $(\ref{eq:1.1})$ 
has at most one possible solution in positive integers. Moreover, 
if $v=k^{2}l$ where $l$ is square-free, then this solution, if it  
exists, is given by 
\begin{displaymath}
x + y^{2}\sqrt{d} = \epsilon^{l}. 
\end{displaymath}
\end{theorem}

Theorem~\ref{thm:B} was proved by making use of the theory of 
linear forms in logarithms of algebraic numbers and some algebraic 
number theory. 

Other special cases of (\ref{eq:1.1}) are known. For $d=2$,  
Ljunggren has shown that (\ref{eq:1.1}) has precisely two 
solutions, $(1,1)$ and $(239,13)$, in positive integers. 
Cohn has also contributed to the study of these equations, 
see \cite[Theorem 7]{Cohn}, for example. 

In \cite{Chen1}, the first author used the method of Thue and Siegel 
to give a new proof that $(1,1)$ and $(239,13)$ are the only solutions 
in positive integers of (\ref{eq:1.1}) for $d=2$. In that paper, 
Pade approximations were used to find an effective improvement of 
Liouville's theorem which was, in turn, used to solve 
$X^{2} + 1 = 2Y^{4}$ via a Thue equation.  

Here we extend the method of \cite{Chen1}. We will establish 
an effective improvement of Liouville's theorem on the measure of 
irrationality for certain algebraic numbers of degree 4, and then 
use such a result to solve a family of Thue equations. In this way, 
we will prove: 

\begin{theorem}
\label{thm:E}
Let $t \geq 1$ be a rational integer with $t \neq 3$. 
Consider the Thue equation 
\begin{equation}
\label{eq:3.1}
P_{t}(X,Y) = X^{4}-tX^{3}Y-6X^{2}Y^{2}+tXY^{3}+Y^{4} = \pm 1. 
\end{equation}

For $t=1$, if $(x,y)$ is an integer solution of $(\ref{eq:3.1})$ then 
$(x,y) \in \{(-2,1)$, $(-1,-2)$, $(-1,0)$, $(0,\pm 1)$, $(1,0)$, $(1,2)$,
$(2,-1) \}$. 

For $t=4$, if $(x,y)$ is an integer solution of $(\ref{eq:3.1})$ then 
$(x,y) \in \{(-3,2)$, $(-2,3)$, $(-1,0)$, $(0,\pm 1)$, $(1,0)$, $(2,3)$,
$(3,-2) \}$. 

For $t = 2$ or $t \geq 5$, if $(x,y)$ is an integer solution 
of $(\ref{eq:3.1})$ then $(x,y) = (\pm 1, 0)$ or $(0,\pm 1)$. 
\end{theorem}

We exclude $t=3$ since $P_{3}(X,Y)=(X^{2}+XY-Y^{2})(X^{2}-4XY-Y^{2})$ 
and hence these Thue equations are quite easy to solve. 

Finally, in Section~4, we use this result to prove the following 
one which concerns (\ref{eq:1.1}). 

\begin{theorem}
\label{thm:C}
Let $(u,v)$ be the fundamental solution of the Pell equation
$X^{2} + 1 = dY^{2}$. If $d \geq 3$, then the equation 
$(\ref{eq:1.1})$ has at most one solution in positive integers. 
If this solution $(x,y)$ exists, we have $v=y^{2}$. 
\end{theorem} 

We note that Lettl and Peth\H{o} \cite{LP} have independently 
proved Theorem~\ref{thm:E} using lower bounds for linear forms 
in two logarithms in the same manner as Thomas \cite{Thom}. 
Actually, they also consider the case when the constant is 
$\pm 4$, but this can be reduced to studying the Thue equations 
above (see Lemma~1 of \cite{LP}). 

As our proof is completely different, based on the explicit 
construction of ``good'' rational approximations to certain 
algebraic numbers by hypergeometric methods, and this appears 
to be the first time that such methods are used to solve a 
family of Thue equations, we feel that there is reason to 
present our own proof here. 

\vspace{3.0mm}

\begin{center}
2. {\sc Preliminaries}
\end{center}

\setcounter{equation}{0}
\setcounter{lemma}{0}
\addtocounter{section}{1}

\vspace{3.0mm}

We start with some notation. 

\vspace{3.0mm}

\noindent
{\bf Notation.} For positive integers $n$ and $r$, we put  
\begin{displaymath}
X_{n,r}(X) = {} _{2}F_{1}(-r,-r-1/n,1-1/n,X),  
\end{displaymath}
where $_{2}F_{1}$ denotes the classical hypergeometric function. 

We use $X_{n,r}^{*}$ to denote the homogeneous 
polynomials derived from these polynomials, so that  
\begin{displaymath}
X_{n,r}^{*}(X,Y) = Y^{r} X_{n,r}(X/Y).  
\end{displaymath}

\vspace{3.0mm}

We now present the following important lemma of Thue.

\begin{lemma}[Thue \cite{Thue}] 
\label{lem:2.3}
Let $P(X)$ be a polynomial of degree $n \geq 2$ and assume 
that there is a quadratic polynomial $U(X)$ with non-zero 
discriminant such that 
\begin{equation}
\label{eq:2.8}
U(X) P''(X) - (n-1)U'(X)P'(X) + \frac{n(n-1)}{2} U''(X)P(X) = 0. 
\end{equation}

We write 
\begin{eqnarray*}
Y_{1}(X) & = & 2U(X)P'(X)-nU'(X)P(X), \\  
h        & = & \frac{n^{2}-1}{4} \left( U'(X)^{2} - 2U(X)U''(X) \right) 
\hspace{5mm} \mbox{ and } \hspace{5mm} 
\lambda = \frac{h}{n^{2}-1}. 
\end{eqnarray*}

Let us define two families of polynomials $A_{r}(X)$ 
and $B_{r}(X)$ by the initial conditions
\begin{eqnarray*}
A_{0}(X) & = & \frac{2h}{3},  \hspace{1.0cm} 
A_{1}(X)   =   \frac{2(n+1)}{3} 
	       \left( U(X)P'(X) - \frac{n-1}{2}U'(X)P(X) \right), \\ 
B_{0}(X) & = & \frac{2hX}{3},  \hspace{1.0cm} 
B_{1}(X)   =   XA_{1}(X) - \frac{2(n+1)U(X)P(X)}{3},   
\end{eqnarray*}
and, for $r \geq 1$, by the recurrence equations 
\begin{align}
\label{eq:2.9}
\lambda (n(r+1)-1) A_{r+1}(X) &= \left( r + \frac{1}{2} \right) Y_{1}(X) A_{r}(X) 
      - (nr+1) P^{2}(X) A_{r-1}(X), \nonumber \\
\lambda(n(r+1)-1)B_{r+1}(X) &= \left( r + \frac{1}{2} \right) Y_{1}(X) B_{r}(X) 
      - (nr+1)P^{2}(X)B_{r-1}(X).  
\end{align}

\noindent
{\rm (i)} For any root $\beta$ of $P(X)$, 
\begin{displaymath}
\beta A_{r}(X) - B_{r}(X) = C_{r}(X) = (X - \beta)^{2r+1} S_{r}(X), 
\end{displaymath}
where $S_{r}(X)$ is a polynomial. 

\noindent
{\rm (ii)} Put 
\begin{displaymath}
z(X) = \frac{1}{2} \left( \frac{Y_{1}(X)}{2n \sqrt{\lambda}} + P(X) \right),  
\hspace{5mm}                          
u(X) = \frac{1}{2} \left( \frac{Y_{1}(X)}{2n \sqrt{\lambda}} - P(X) \right)  
\hspace{3.0mm} \mbox{ and } \hspace{3.0mm}                         
w(X) = \frac{z(X)}{u(X)}.
\end{displaymath}

Then 
\begin{eqnarray*}
(\sqrt{\lambda})^{r}A_{r}(X) & = & a(X)X_{n,r}^{*}(z,u) - b(X)X_{n,r}^{*}(u,z) 
\mbox{ and } \\
(\sqrt{\lambda})^{r}B_{r}(X) & = & c(X)X_{n,r}^{*}(z,u) - d(X)X_{n,r}^{*}(u,z), 
\end{eqnarray*}
where
\begin{eqnarray*}
a(X) & = & \left( \frac{(n-1)\sqrt{\lambda}}{2P(X)} \right) A_{1}(X)  
	   - \left( \frac{Y_{1}(X)}{4\sqrt{\lambda}P(X)} - \frac{1}{2} \right) 
	     A_{0}(X), \\ 
b(X) & = & \left( \frac{(n-1)\sqrt{\lambda}}{2P(X)} \right) A_{1}(X)  
	   - \left( \frac{Y_{1}(X)}{4\sqrt{\lambda}P(X)} + \frac{1}{2} \right) 
	     A_{0}(X), \\ 
c(X) & = & \left( \frac{(n-1)\sqrt{\lambda}}{2P(X)} \right) B_{1}(X)  
	   - \left( \frac{Y_{1}(X)}{4\sqrt{\lambda}P(X)} - \frac{1}{2} \right) 
	     B_{0}(X) \mbox{ and } \\
d(X) & = & \left( \frac{(n-1)\sqrt{\lambda}}{2P(X)} \right) B_{1}(X)  
	   - \left( \frac{Y_{1}(X)}{4\sqrt{\lambda}P(X)} + \frac{1}{2} \right) 
	     B_{0}(X). \\ 
\end{eqnarray*}
\end{lemma}

These results can be found in Thue \cite[Theorem and equations 
35--47]{Thue} or Chudnovsky \cite{Chud} (see, in particular, 
Lemma~7.1 and the remarks that follow (pages 364--366)). 

We have added two extra hypotheses, requiring that the degree 
of $P(X)$ be at least two and that the discriminant of $U(X)$ be 
non-zero. Since $h$ is equal to $(n^{2}-1)/4$ times the discriminant 
of $U(X)$ and since $\lambda$ is a multiple of $h$, we need these 
conditions to ensure that we do not fall into degenerate cases with 
the $A_{i}$'s and $B_{i}$'s. 

Also notice that there are some differences in notation between the 
lemma above, which is similar to Chudnovsky's \cite{Chud}, and that 
of Thue. In particular, here, Thue's $\alpha$ and $F$ are replaced 
by $\beta$ and $P$, $n$ and $r$ are switched from \cite{Thue}, our 
$A_{i}$'s and $B_{i}$'s are $2(n-1)/3$ times Thue's and we label 
Thue's $U_{r}(z,y)$ as $X_{n,r}^{*}(z,u)$. Also, what we call 
$Y_{1}(X)$ and $u(X)$ respectively, correspond to $2H(x)$ and 
$y(x)$ respectively in Thue's paper. 

We now give two lemmas which will give us the remainder of 
our approximations in a nice form that can be easily bounded 
from above. 

\begin{lemma}
\label{lem:2.4}
Let $m$ and $n$ be positive integers and suppose 
$\alpha$ is a real number. Define 
\begin{eqnarray*}
p_{m}(X) = \sum_{\nu=0}^{m} {m-\alpha \choose m-\nu} 
			    {n+\alpha \choose \nu} X^{\nu} 
\hspace{3.0mm} \mbox{ and } \hspace{3.0mm}                                  
q_{n}(X) = \sum_{\nu=0}^{n} {m-\alpha \choose \nu} 
			    {n+\alpha \choose n-\nu} X^{\nu}.  
\end{eqnarray*}

Given a complex number $x$, we let $C$ denote the straight line from 
$1$ to $x$. If $0$ is not on $C$, then 
\begin{displaymath}
x^{\alpha} q_{n}(x) - p_{m}(x) 
= \alpha {m-\alpha \choose m} {n+\alpha \choose n}        
  \int_{C} (t-x)^{m} (1-t)^{n} t^{\alpha-m-1} dt. 
\end{displaymath}
\end{lemma}

\begin{proof}
Put 
\begin{displaymath}
r(x) = x^{\alpha} q_{n}(x) - p_{m}(x).  
\end{displaymath}

It is a routine matter to verify that, when $k=0,1,\ldots,m$, 
\begin{displaymath}
\frac{r^{(k)}(x)}{k!}
= \sum_{\nu=0}^{n} {m-\alpha \choose \nu} {n+\alpha \choose n-\nu} 
		   {\nu+\alpha \choose k} x^{\alpha+\nu-k}             
  - \sum_{\nu=k}^{m} {m-\alpha \choose m-\nu} {n+\alpha \choose \nu} 
		   {\nu \choose k} x^{\nu-k}. 
\end{displaymath}

Note that 
\begin{displaymath}
{n+\alpha \choose n-\nu} {\nu+\alpha \choose k} 
= {n-k+\alpha \choose n-\nu} {n+\alpha \choose k}.  
\end{displaymath}

Thus, using the famous Vandermonde formula 
\begin{displaymath}
\sum_{\nu=0}^{k} {x \choose \nu} {y \choose k-\nu} 
= {x+y \choose k},                               
\end{displaymath}
where $x$ and $y$ are real numbers, we get 
\begin{eqnarray*}
\sum_{\nu=0}^{n} {m-\alpha \choose \nu} {n+\alpha \choose n-\nu} 
		 {\nu+\alpha \choose k} 
& = & {n+\alpha \choose k} 
      \sum_{\nu=0}^{n} {m-\alpha \choose \nu} {n-k+\alpha \choose n-\nu} \\ 
& = & {n+\alpha \choose k} {m+n-k \choose n}.  
\end{eqnarray*}

We also have 
\begin{displaymath}
{n+\alpha \choose \nu+k} {\nu+k \choose k} 
= {n+\alpha \choose k} {n-k+\alpha \choose \nu}.  
\end{displaymath}

Applying this identity along with Vandermonde's formula again, we get 
\begin{eqnarray*}
\sum_{\nu=k}^{m} {m-\alpha \choose m-\nu} {n+\alpha \choose \nu} 
		 {\nu \choose k} 
& = & {n+\alpha \choose k} 
      \sum_{\nu=0}^{m-k} {m-\alpha \choose m-k-\nu} 
			 {n-k+\alpha \choose \nu} \\ 
& = & {n+\alpha \choose k} {m+n-k \choose n}. 
\end{eqnarray*}

Hence $r^{(k)}(1)/k! = 0$ for $k=0,1,\ldots,m$. 

Further computation shows 
\begin{eqnarray*}
{m-\alpha \choose \nu} {n+\alpha \choose n-\nu} 
{\nu+\alpha \choose m+1} 
& = & {n \choose \nu} 
      \frac{(m-\alpha) \cdots (m-\alpha-\nu+1)
	    (n+\alpha) \cdots (\alpha+\nu-m)}{(m+1)! \, n!}  \\ 
& = & (-1)^{\nu-m} \frac{\alpha}{m+1} {m-\alpha \choose m} 
      {n+\alpha \choose n} {n \choose \nu}.  
\end{eqnarray*}

Thus 
\begin{eqnarray*}
\frac{r^{(m+1)}(x)}{(m+1)!}
& = & \sum_{\nu=0}^{n} {m-\alpha \choose \nu} {n+\alpha \choose n-\nu} 
		   {\nu+\alpha \choose m+1} x^{\alpha+\nu-m-1} \\             
& = & (-1)^{m} \frac{\alpha}{m+1} {m-\alpha \choose m} 
      {n+\alpha \choose n} x^{\alpha-m-1} (1-x)^{n}. 
\end{eqnarray*}

Expanding $r(x)$ into its Taylor series with remainder centred around 
$x_{0}=1$, we have 
\begin{displaymath}
r(x) = r(1)+r'(1)(x-1) + \ldots + \frac{1}{m!}r^{(m)}(1)(x-1)^{m} 
       + \frac{1}{m!} \int_{1}^{x} r^{(m+1)}(t)(x-t)^{m} dt.  
\end{displaymath}

Hence 
\begin{displaymath}
r(x) = \alpha {m-\alpha \choose m} {n+\alpha \choose n}        
       \int_{C} t^{\alpha-m-1}(t-x)^{m}(1-t)^{n} dt,   
\end{displaymath}
and the lemma follows. 
\end{proof}

\begin{lemma}
\label{lem:2.5} 
Let $n \geq 2$ and $r$ be positive integers and let $\beta$, $\lambda$, $a(X)$,
$b(X)$, $c(X)$, $C_{r}(X)$, $d(X),u(X)$, $w(X),X_{n,r}^{*}(X,Y)$ and $z(X)$ be
as in Lemma~$\ref{lem:2.3}$. Put 
\begin{displaymath}
R_{n,r}(x) = \frac{\Gamma(r+1+1/n)}{r! \Gamma(1/n)} 
	     \int_{1}^{x} (1-t)^{r}(t-x)^{r}t^{-(r+1-1/n)} dt,   
\end{displaymath}             
where the integration path is the straight line from $1$ to $x$. 

If $x$ is a complex number such that $w(x)$ is not a negative number or zero,
then
\begin{eqnarray*}
({\sqrt{\lambda}})^{r}C_{r}(x) 
& = & \left( \beta \left( a(x)w(x)^{1/n} - b(x) \right) 
	     - \left( c(x)w(x)^{1/n} - d(x) \right) \right) 
      X_{n,r}^{*}(u(x),z(x)) \\ 
&   & - \left( \beta a(x) - c(x) \right) u(x)^{r}R_{n,r}(w(x)).  
\end{eqnarray*}
\end{lemma}

\begin{proof}
Letting $\alpha=1/n$, it is easy to verify that  
\begin{displaymath}
p_{r}(x) = {r - 1/n \choose r} X_{n,r}(x) 
\hspace{3.0mm} \mbox{ and } \hspace{3.0mm} 
q_{r}(x) = {r - 1/n \choose r} x^{r} X_{n,r}(1/x). 
\end{displaymath}

Observe that here, and in what follows, $X_{n,r}(x)$ is a polynomial of
degree $r$, so $x^{r} X_{n,r}(1/x)$ is also a polynomial of degree $r$ and
hence well-defined.

Thus from Lemma~\ref{lem:2.4} and the definition of $R_{n,r}(x)$, 
we have 
\begin{displaymath}
x^{1/n}x^{r}X_{n,r}(1/x)-X_{n,r}(x) 
= \frac{\Gamma(r+1+1/n)}{r! \, \Gamma(1/n)} 
  \int_{1}^{x} (1-t)^{r}(t-x)^{r}t^{1/n-r-1} dt 
= R_{n,r}(x). 
\end{displaymath}
  
Substituting $x=w$ into this expression, we have 
$w^{1/n}X_{n,r}^{*}(u,z) = X_{n,r}^{*}(z,u) + u^{r}R_{n,r}(w)$. 
So replacing $X_{n,r}^{*}(z,u)$ by 
$w^{1/n}X_{n,r}^{*}(u,z) - u^{r}R_{n,r}(w)$ 
in the expressions for $(\sqrt{\lambda})^{r}A_{r}(X)$ and 
$(\sqrt{\lambda})^{r}B_{r}(X)$ in Lemma~\ref{lem:2.3}(ii) 
and then applying the relation in Lemma~\ref{lem:2.3}(i), 
our result follows. 
\end{proof}

\vspace{3.0mm}

We now give a generalization of a result of Baker \cite[Lemma 3]{Baker},  
itself an improvement of a result of Siegel \cite[Hilfssatz 5]{Siegel}, 
in a form which is suitable for our needs here. 

\begin{lemma} 
\label{lem:2.6}
Suppose $j = \pm 1$ and $n$ and $r$ are positive integers. We define 
\begin{displaymath}
\mu_{n} = \prod_{\stackrel{\displaystyle p|n}{p, {\rm prime}}} 
p^{1/(p-1)}.  
\end{displaymath}

Then the coefficients of the polynomial 
\begin{displaymath}
{2r \choose r} {} _{2}F_{1}(-r,-r+j/n, -2r, n \mu_{n}X) 
\end{displaymath}
are algebraic integers. 
\end{lemma}

\begin{proof}
Let 
\begin{displaymath}
p(X) = {2r \choose r} {} _{2}F_{1}(-r,-r+j/n, -2r, n \mu_{n}X).  
\end{displaymath}

By definition of the hypergeometric function, we have 
\begin{displaymath}
p(X) = \sum_{s=0}^{r} \frac{l_{s}}{s!} n^{-s} {2r-s \choose r} 
       (-n \mu_{n}X)^{s}   
     = \sum_{s=0}^{r} \frac{l_{s}}{s!} \mu_{n}^{s} {2r-s \choose r} 
       (-X)^{s}     
\end{displaymath}
where  
\begin{displaymath}
l_{s} = \prod_{k=r-s+1}^{r} (kn-j). 
\end{displaymath}

Defining 
\begin{displaymath}
\sigma_{s} = \prod_{\stackrel{\displaystyle p|n}{p, {\rm prime}}} 
p^{[s/(p-1)]},   
\end{displaymath}
for non-negative $s$, from Lemma 4.1 of \cite{Chud}, we see 
that $l_{s} \sigma_{s}/s!$ is a rational integer. Since 
\begin{displaymath}
\frac{\mu_{n}^{s}}{\sigma_{s}} 
= \prod_{\stackrel{\displaystyle p|n}{p, {\rm prime}}} p^{s/(p-1)-[s/(p-1)]}  
\end{displaymath}
and $s/(p-1)-[s/(p-1)]$ is a non-negative rational number, 
$\mu_{n}^{s}/\sigma_{s}$ is an algebraic integer. Therefore, 
\begin{displaymath}
\frac{l_{s}}{s!}\sigma_{s} \frac{\mu_{n}^{s}}{\sigma_{s}} 
{2r-s \choose r}
= \frac{l_{s}}{s!} \mu_{n}^{s} {2r-s \choose r} 
\end{displaymath}
is an algebraic integer and the lemma follows. 
\end{proof}

\begin{lemma}
\label{lem:2.1}
Let $w = e^{i \varphi}, 0 < \varphi < \pi$ and put 
$\sqrt{w}=e^{i \varphi/2}$. Let $n$ and $r$ be positive integers. 
Define $R_{n,r}(x)$ as in Lemma~$\ref{lem:2.5}$. Then 
\begin{displaymath}
\left| R_{n,r}(w) \right| 
\leq \frac{\Gamma(r+1+1/n)}{r! \, \Gamma(1/n)}
     \varphi { \left| 1 - \sqrt{w} \right| }^{2r}. 
\end{displaymath}
\end{lemma}

\begin{proof}
By Cauchy's theorem, 
\begin{displaymath}
R_{n,r}(w) 
= \frac{\Gamma(r+1+1/n)}{r! \, \Gamma(1/n)}
  \int_{C} { \left( (1-t)(t-w) \right) }^{r} t^{1/n-r-1} dt,   
\end{displaymath}
where 
\begin{displaymath}
C = \{ t \, | \, t=e^{i \theta}, 0 \leq \theta \leq \varphi \}. 
\end{displaymath}

Put 
\begin{displaymath}
f(t) = \frac{(1-t)(t-w)}{t}  
\hspace{5 mm} \mbox{ and } \hspace{5 mm}
g(t) = t^{1/n-1}. 
\end{displaymath}

Define 
\begin{displaymath}
F(\theta) = { \left| f \left( e^{i \theta} \right) \right| }^{2}, 
\end{displaymath}
so 
\begin{displaymath}
F(\theta) = 4 (1 - \cos \theta)(1 - \cos (\theta - \varphi))  
\hspace{5.0mm} \mbox{ for $0 \leq \theta \leq \varphi$.}  
\end{displaymath}

A simple calculation shows that 
\begin{displaymath}
F'(\theta) = -16 \sin \left( \theta - \frac{\varphi}{2} \right) 
	     \sin \left( \frac{\theta}{2} \right)
	     \sin \left( \frac{\varphi - \theta}{2} \right).   
\end{displaymath}

The only values of $0 \leq \theta \leq \varphi$ with $F'(\theta)=0$ 
are $\theta=0,\varphi/2$ and $\varphi$. It is easy to check that 
\begin{displaymath}
F(\theta) \leq F(\varphi/2) 
= 4 { \left(  1 - \cos \frac{\varphi}{2} \right) }^{2} 
= { \left| 1 - \sqrt{w} \right| }^{4}, 
\end{displaymath}
and hence                                             
\begin{displaymath}
\left| \int_{C} f(t)^{r} g(t) dt \right|
\leq \int_{0}^{\varphi} 
     { \left| f \left( e^{i \theta} \right) \right| }^{r} 
     \left| g \left( e^{i \theta} \right) \right| d \theta 
\leq \varphi { \left| 1 - \sqrt{w} \right| }^{2r}. 
\end{displaymath}

The lemma follows. 
\end{proof}

\begin{lemma}
\label{lem:2.8}
Let $u$ and $z$ be complex numbers with $w = z/u = e^{i \varphi}$ 
where $0 < \varphi < \pi$. Then 
\begin{eqnarray*}
\left| X_{n,r}^{*}(z,u) \right| 
& \leq & 4 |u|^{r} \frac{\Gamma(1-1/n) \, r!}{\Gamma(r+1-1/n)} 
	 { \left| 1 + \sqrt{w} \right| }^{2r-2} \mbox{ and } \\ 
\left| X_{n,r}^{*}(u,z) \right| 
& \leq & 4 |z|^{r} \frac{\Gamma(1-1/n) \, r!}{\Gamma(r+1-1/n)} 
	 { \left| 1 + \sqrt{w} \right| }^{2r-2}. 
\end{eqnarray*}
\end{lemma}

\begin{proof}
Recall that $X_{n,r}^{*}(z,u) = u^{r} {} _{2}F_{1}(-r,-r-1/n,1-1/n,w)$.  

We noted in the proof of Lemma~\ref{lem:2.5} that 
\begin{displaymath}
\sum_{k=0}^{r} {r-1/n \choose r-k} {r+1/n \choose k} w^{k} 
= {r-1/n \choose r} {} _{2}F_{1}(-r,-r-1/n,1-1/n,w).  
\end{displaymath}

So, by the binomial theorem and Cauchy's residue theorem, 
\begin{displaymath}
_{2}F_{1} (-r,-r-1/n,1-1/n,w) 
= (-1)^{r} \frac{\Gamma(1-1/n) \, r!}{2 \pi i \Gamma(r+1-1/n)}
  \int_{C} t^{-r-1} (1-t)^{r-1/n} (1-wt)^{r+1/n} dt, 
\end{displaymath}  
where $C$ is any path which encircles the origin 
once in the positive sense. 

We now focus on bounding the absolute value of this integral 
from above. 

Put 
\begin{displaymath}
f(t) = \frac{(1-t)(1-wt)}{t} 
\hspace{3.0mm} \mbox{ and } \hspace{3.0mm}
g(t) = \frac{(1-t)^{1-1/n}(1-wt)^{1+1/n}}{t^{2}}. 
\end{displaymath}

Now
\begin{displaymath}
F(\theta) = { \left| f \left( e^{i \theta} \right) \right| }^{2} 
	  = 4 (1 - \cos \theta)(1 - \cos (\theta + \varphi))  
\hspace*{5mm} \mbox{ for $0 \leq \theta < 2 \pi$.}
\end{displaymath}

We have 
\begin{displaymath}
F'(\theta) = -16 \sin \left( \theta + \frac{\varphi}{2} \right) 
	     \sin \left( \frac{\theta}{2} \right)  
	     \sin \left( \frac{- \varphi - \theta}{2} \right),  
\end{displaymath}
so that the only values of $0 \leq \theta < 2 \pi$ for which 
$F'(\theta)=0$ are $\theta = 0, \pi - \varphi/2, 2 \pi - \varphi$ 
and $2 \pi - \varphi/2$. It is easy to check that if 
$0 < \varphi < \pi$ then 
$F(\theta) \leq F(\pi - \varphi/2) = |1+\sqrt{w}|^{4}$. 

Note that 
\begin{displaymath}
\left| g \left( e^{i\theta} \right) \right| \leq 4. 
\end{displaymath}

Thus, by Cauchy's theorem and the Cauchy-Schwartz inequality, we get  
\begin{displaymath}
\left| \int_{C} t^{-r-1} (1-t)^{r-1/n} (1-wt)^{r+1/n} dt \right| 
\leq \int_{0}^{2\pi} F(\theta)^{(r-1)/2} |g(e^{i\theta})| d \theta
\leq 8 \pi { \left| 1+\sqrt{w} \right| }^{2(r-1)}.  
\end{displaymath}

Therefore, 
\begin{displaymath}
\left| {} _{2}F_{1} (-r,-r-1/n,1-1/n,w) \right| 
\leq \frac{\Gamma(1-1/n) \, r!}{2 \pi \Gamma(r+1-1/n)}
     8 \pi \left| 1 + \sqrt{w} \right|^{2r-2}.
\end{displaymath}  

To prove that the second inequality holds, we note that since $w=z/u$ is on the
unit circle, we have that $|u|=|z|$ and $w^{-1}=\overline{w}$. Therefore,
${} _{2}F_{1} \left( -r,-r-1/n,1-1/n,w^{-1} \right)$ is the complex conjugate
of ${} _{2}F_{1} (-r,-r-1/n,1-1/n,w)$. The validity of the first inequality now
implies the second.  
\end{proof}

\begin{lemma}
\label{lem:2.7} 
Let $A_{r}(X),B_{r}(X),P(X)$ and $U(X)$ be defined as in 
Lemma~$\ref{lem:2.3}$ and let $a,b,c$ and $d$ be complex 
numbers satisfying $ad-bc \neq 0$. Define 
\begin{displaymath}
K_{r}(X) = a A_{r}(X) + b B_{r}(X) 
\mbox{\hspace{5.0mm} and \hspace{5.0mm}} 
L_{r}(X) = c A_{r}(X) + d B_{r}(X).  
\end{displaymath}

If $U(x)P(x) \neq 0$ for a given $x$, then 
\begin{displaymath}
K_{r+1}(x)L_{r}(x) \neq K_{r}(x)L_{r+1}(x),  
\end{displaymath}
for all $r \geq 0$.
\end{lemma}

\begin{proof}
We prove this lemma by induction on $r$. 

Using the relations between the $K_{i}$'s and $L_{i}$'s and 
the $A_{i}$'s and $B_{i}$'s, as well as the definitions of 
$A_{0}(X),A_{1}(X),B_{0}(X)$ and $B_{1}(X)$ in Lemma~\ref{lem:2.3}, 
we obtain 
\begin{displaymath}
K_{1}(x) L_{0}(x) - K_{0}(x) L_{1}(x) 
= \frac{4h(n+1)P(x)U(x)(ad-bc)}{9}, 
\end{displaymath}
where $h$ is as in Lemma~\ref{lem:2.3}. 

Recalling that our assumptions on $P(X)$ and $U(X)$ in 
Lemma~\ref{lem:2.3} ensure that $h$ and $\lambda$ are 
non-zero, so is $K_{1}(x) L_{0}(x) - K_{0}(x) L_{1}(x)$.  

Now assume that the lemma holds for all non-negative integers 
up to $r-1$. 

Since $A_{r}(X)$ and $B_{r}(X)$ satisfy the recurrence (\ref{eq:2.9}), 
the sequences of numbers $\left( K_{r}(x) \right)$ and 
$\left( L_{r}(x) \right)$ satisfy the same recurrence, namely 
\begin{eqnarray*}
\lambda (n(r+1)-1) K_{r+1}(x) 
& = & \left( r + \frac{1}{2} \right) Y_{1}(x) K_{r}(x) 
      - (nr+1) P^{2}(x) K_{r-1}(x) \mbox{ and } \\ 
\lambda (n(r+1)-1) L_{r+1}(x) 
& = & \left( r + \frac{1}{2} \right) Y_{1}(x) L_{r}(x) 
      - (nr+1) P^{2}(x) L_{r-1}(x),  
\end{eqnarray*}
where $\lambda$ is as in Lemma~\ref{lem:2.3} and $n$ is the 
degree of $P(X)$.

Multiplying both sides of the first equation by $L_{r}(x)$,  
both sides of the second equation by $K_{r}(x)$ and 
subtracting the resulting expressions, we obtain 
\begin{displaymath}
K_{r+1}(x) L_{r}(x) - K_{r}(x) L_{r+1}(x) 
= \frac{(nr+1)P^{2}(x)}{\lambda (n(r+1)-1)} 
\left( K_{r}(x) L_{r-1}(x) - K_{r-1}(x) L_{r}(x) \right).  
\end{displaymath}

The lemma now holds for $r$ by our inductive hypothesis.
\end{proof}

\begin{lemma}
\label{lem:2.9}
Let $\theta \in \bbR$. Suppose that there exist $k_{0},l_{0} > 0$ 
and $E,Q > 1$ such that for all $r \in \bbN$, there are rational 
integers $p_{r}$ and $q_{r}$ with $|q_{r}| < k_{0}Q^{r}$ and 
$|q_{r} \theta - p_{r}| \leq l_{0}E^{-r}$ satisfying 
$p_{r}q_{r+1} \neq p_{r+1}q_{r}$. Then for any rational 
integers $p$ and $q$ with $|q| \geq 1/(2l_{0})$, we have 
\begin{displaymath}
\left| \theta - \frac{p}{q} \right| > \frac{1}{c |q|^{\kappa+1}}, 
\mbox{ where $c=2k_{0}Q(2l_{0}E)^{\kappa}$ and 
	     $\kappa = \displaystyle \frac{\log Q}{\log E}$.} 
\end{displaymath}             
\end{lemma}

\begin{proof}
Let $p,q \in \bbZ$ with $|q| \geq 1/(2l_{0}) > 0$ be given. 
Choose 
$\displaystyle n_{0} = \left[ \frac{\log(2l_{0}|q|)}{\log E} \right]+1$.  
We have $n_{0} \geq 1$ since $2l_{0}|q| \geq 1$. From 
$\log (2l_{0}|q|)/\log E < n_{0}$, one deduces that 
$l_{0}E^{-n} \leq l_{0}E^{-n_{0}} < 1/(2|q|)$ for $n \geq n_{0}$. 

If, for some $n \geq n_{0}$, we have $q_{n}=0$, then 
$|q_{n} \theta - p_{n}| \leq l_{0}E^{-n} < 1/(2|q|) < 1$ 
would imply that $p_{n}=0$, contradicting the supposition 
that $p_{n}q_{n+1} \neq p_{n+1}q_{n}$. Therefore, $q_{n} \neq 0$ 
for all $n \geq n_{0}$. 

Notice that for any $n \geq n_{0}$ with $p/q \neq p_{n}/q_{n}$, 
we have 
\begin{displaymath}
\left| \theta - \frac{p}{q} \right| 
\geq \left| \frac{p_{n}}{q_{n}} - \frac{p}{q} \right| 
     - \left| \theta - \frac{p_{n}}{q_{n}} \right| 
\geq \frac{1}{|qq_{n}|} - \frac{l_{0}}{E^{n}|q_{n}|}     
> \frac{1}{2|qq_{n}|}. 
\end{displaymath}

The choice of $n_{0}$ yields $Q^{n_{0}} \leq (2El_{0}|q|)^{\kappa}$. 
If $p/q \neq p_{n_{0}}/q_{n_{0}}$ then we can conclude, as in the 
proof of Lemma~2.9 of \cite{Chen1}, and obtain 
\begin{displaymath}
\left| \theta - \frac{p}{q} \right| 
> \frac{1}{2|q|k_{0}Q^{n_{0}}} 
\geq \frac{1}{2k_{0}(2El_{0})^{\kappa}|q|^{\kappa+1}}.  
\end{displaymath}

If $p/q = p_{n_{0}}/q_{n_{0}}$ then we have 
$p/q \neq p_{n_{0}+1}/q_{n_{0}+1}$ and obtain 
\begin{displaymath}
\left| \theta - \frac{p}{q} \right| 
> \frac{1}{2|q|k_{0}Q^{n_{0}+1}} 
\geq \frac{1}{2k_{0}Q(2El_{0})^{\kappa}|q|^{\kappa+1}}.  
\end{displaymath}
\end{proof}

\vspace{3.0mm}

\noindent
{\bf Note.} This lemma is very much like Lemma~3 of \cite{Chen1},  
except that the lower bound for $|q|$ is weaker here and we have 
included an extra factor of $Q$ in $c$. This extra factor seems 
required by the consideration in the proof that $p/q=p_{n}/q_{n}$ 
is possible for some $n$, a consideration omitted from the proof 
of Lemma~3 in \cite{Chen1}. It was the referee who spotted this 
omission and we are thankful for his/her attentive reading in 
observing this and suggesting a way of rectifying it. 

\begin{lemma}
\label{lem:2.11}
Let $\alpha$ and $\beta$ be algebraic integers such that 
$\alpha /\beta$ is not a root of unity and $\alpha+\beta$ and 
$\alpha \beta$ are relatively prime rational integers. Define 
the sequence ${ \left( V_{m} \right) }_{m=-\infty}^{\infty}$ by 
\begin{displaymath}
V_{m} = \frac{\alpha^{m}-\beta^{m}}{\alpha - \beta}. 
\end{displaymath}

\noindent
{\rm (i)} If $\alpha$ and $\beta$ are roots of the polynomial 
$X^{2}-aX-b$ then $V_{m}=aV_{m-1}+bV_{m-2}$ for $m \geq 2$. 

\noindent
{\rm (ii)} If $m$ and $n$ are positive integers with $d=(m,n)$ 
then $(V_{m},V_{n})=V_{d}$. 

\noindent
{\rm (iii)} If the minimal polynomial of $\alpha$ and $\beta$ 
over $\bbZ$ is $X^{2}-aX-1$ where $a$ is a positive even integer 
then, for any relatively prime odd positive integers $m$ and $n$, 
\begin{displaymath}
\left( \frac{V_{m}}{V_{n}} \right) = 1, 
\end{displaymath}
where $(\cdot / \cdot)$ denotes the Jacobi symbol. 
\end{lemma}

\begin{proof}
(i) This is a well-known fact from the theory of linear 
recurrence sequences (see, for example, Chapter C of \cite{ST}). 

(ii) This follows using the ideas in the proof of Theorem~179(iv) 
of \cite{HW}. 

(iii) We first need to define the associated sequence 
${ \left( U_{m} \right) }_{m=-\infty}^{\infty}$ by 
\begin{displaymath}
U_{m} = \alpha^{m}+\beta^{m},   
\end{displaymath}
for all integers $m$. 

Now we take a moment to establish some simple relationships 
that elements of these sequences satisfy. 

Notice that $(-1)^{m} U_{-m} = (\alpha \beta)^{m} U_{-m} = U_{m}$ 
and $(-1)^{m} V_{-m} = -V_{m}$.  

We can quickly verify from the definitions of the $U_{i}$'s 
and $V_{i}$'s that 
\begin{eqnarray}
\label{eq:sum}
2V_{m+n} & = & U_{m}V_{n}+U_{n}V_{m}, \nonumber \\
U_{m}^{2}-(a^{2}+4)V_{m}^{2} & = & 4(-1)^{m} \hspace{3.0mm} \mbox{ and }  \nonumber \\
U_{2km} - 2(-1)^{km} & = & (a^{2}+4)V_{km}^{2}, 
\end{eqnarray}
upon noting that $(\alpha-\beta)^{2}
=(\alpha+\beta)^{2}-4\alpha\beta=a^{2}+4$. 

Since $V_{m}$ divides $V_{km}$ and $a$ is even, 
from the last formula we see that 
$U_{2km} \equiv 2(-1)^{km} \bmod 2V_{m}$, so 
the $U_{i}$'s are even and 
\begin{equation}
\label{eq:acong}
\frac{U_{2km}}{2} \equiv (-1)^{km} \bmod V_{m}. 
\end{equation}

By (\ref{eq:sum}), we have $2V_{2n}=2U_{n}V_{n}$. Since 
$U_{n}$ is even, so is $V_{2n}$. This, together with 
the fact that $V_{1}=1$ and the recurrence relation 
$V_{m+1}=aV_{m}+V_{m-1}$ with $a$ even, implies that 
$V_{2n+1} \equiv 1 \bmod 4$. 

We obtain $2V_{2km+r}=U_{r}V_{2km}+U_{2km}V_{r}$ from (\ref{eq:sum}) 
which, together with the observation that $V_{m}$ divides $V_{2km}$, 
implies that $V_{2km+r} \equiv (U_{2km}/2)V_{r} \bmod V_{m}$. 
Combining this congruence with (\ref{eq:acong}), we have 
\begin{equation}
\label{eq:bcongp}
V_{2km+r} \equiv (-1)^{km} V_{r} \bmod V_{m}. 
\end{equation}

Similarly, 
\begin{equation}
\label{eq:bcongm}
V_{2km-r} \equiv (-1)^{km} V_{-r} \equiv (-1)^{km+r-1} V_{r} \bmod V_{m}. 
\end{equation}

We are now ready to prove the assertion in the lemma.

Since $\alpha$ and $\beta$ are real numbers, the $V_{i}$'s 
are positive for positive $i$. And, by part (ii), if $m$ and 
$n$ are relatively prime then so are $V_{m}$ and $V_{n}$, 
since $V_{1}=1$, therefore $(V_{m}/V_{n})$ is defined. 

To prove our claim that 
$$
%\begin{equation}
%\label{eq:jacobi}
\left( \frac{V_{m}}{V_{n}} \right) = 1, 
%\end{equation}
$$
for relatively prime odd positive integers $m$ and $n$, 
we use induction.  

The result is certainly true for $m=1$ as $V_{1}=1$. 

Now suppose that the result holds for all positive odd 
integers $r<m$.  

Since $m$ and $n$ are odd, $V_{m} \equiv V_{n} \equiv 1 \bmod 4$, 
so, by the law of quadratic reciprocity we have                 
$(V_{m}/V_{n})=(V_{n}/V_{m})$. Therefore, without 
loss of generality, we may assume that $m > n$. 

This implies that we can write $m=2kn \pm r$ for some positive 
integer $k$ and an odd positive integer $0 < r < n$. Using 
(\ref{eq:bcongp}), (\ref{eq:bcongm}) and the fact that 
$V_{n} \equiv 1 \bmod 4$, we can show that 
\begin{displaymath}
\left( \frac{V_{m}}{V_{n}} \right) 
= \left( \frac{\pm 1}{V_{n}} \right) 
  \left( \frac{V_{r}}{V_{n}} \right) 
= \left( \frac{V_{r}}{V_{n}} \right).  
\end{displaymath}

Our inductive hypothesis now shows that $(V_{m}/V_{n})=1$.
\end{proof}

\begin{lemma}
\label{lem:2.12}
Let $(u,v)$ be the fundamental solution of the Pell equation 
$X^{2} + 1 = dY^{2}$. If $(\ref{eq:1.1})$ is solvable then 
$v$ is a perfect square.  
\end{lemma}

\begin{proof}
Suppose $(x,y)$ is a solution of (\ref{eq:1.1}), then we have 
\begin{equation}
\label{eq:chena}
y^{2} = \frac{\epsilon^{t}-\overline{\epsilon}^{t}}{2 \sqrt{d}}, 
\end{equation}
where $t \equiv 1 \bmod 2, \epsilon = u + v\sqrt{d}$
and $\overline{\epsilon} = u - v \sqrt{d}$. 

We may assume that $t$ is the least positive integer such that 
(\ref{eq:chena}) is valid. If $t=1$ then the lemma follows so 
we may suppose that $t > 1$ and write $t=ps$ where $p$ is some 
odd prime. From (\ref{eq:chena}), we now obtain  
\begin{equation}
\label{eq:chenb}
y^{2} = \frac{\epsilon^{ps}-\overline{\epsilon}^{ps}}
	     {\epsilon^{s}-\overline{\epsilon}^{s}} 
	\frac{\epsilon^{s}-\overline{\epsilon}^{s}}{2 \sqrt{d}}. 
\end{equation}

Let us define 
\begin{displaymath}
f = \left( \frac{\epsilon^{ps}-\overline{\epsilon}^{ps}} 
		{\epsilon^{s}-\overline{\epsilon}^{s}},
	   \frac{\epsilon^{s}-\overline{\epsilon}^{s}}{2 \sqrt{d}} \right). 
\end{displaymath}
	   
It is well-known that either $f=1$ or $f=p$ (see Lemma~2 of 
\cite{Stewart}, for example). 

If $f=1$ then, from (\ref{eq:chenb}), we see that 
$(\epsilon^{s}-\overline{\epsilon}^{s})/(2 \sqrt{d})$ is 
a perfect square. But since $s<t$ we get a contradiction, 
so $f=p$ and 
$$
%\begin{equation}
%\label{eq:chenc}
py_{1}^{2} = \frac{\epsilon^{ps}-\overline{\epsilon}^{ps}} 
		  {\epsilon^{s}-\overline{\epsilon}^{s}},
%\end{equation}
$$
for some positive integer $y_{1}$.           

To prove that this is impossible, we consider certain 
binary recurrence sequences. There exist integers $g$ 
and $h$ such that 
\begin{displaymath}
\alpha = \epsilon^{s}=h+g\sqrt{d} 
\hspace{3.0mm} \mbox{ and } \hspace{3.0mm} 
\beta  = \overline{\epsilon}^{s}=h-g\sqrt{d}. 
\end{displaymath}

We define the sequence ${ \left( V_{m} \right) }_{-\infty}^{\infty}$ 
by 
\begin{displaymath}
V_{m}  = \frac{\alpha^{m}-\beta^{m}}{\alpha-\beta} 
\end{displaymath}
for all integers $m$. 

We want to show that $pV_{p}$ is not a perfect square, 
so let us assume the opposite. 

Using the definition of the $V_{m}$'s, 
\begin{displaymath}
V_{p} = \sum_{n=0}^{(p-1)/2} {p \choose 2n}h^{2n}(g\sqrt{d})^{p-2n-1},   
\end{displaymath}        
since $p$ is odd. 

For $1 \leq n < p/2$, $p$ divides ${p \choose 2n}$, so 
\begin{displaymath}
V_{p} \equiv (dg^{2})^{(p-1)/2} \bmod p. 
\end{displaymath}

But we are assuming that $pV_{p}$ is a square, $p$ divides 
$V_{p}$ and $dg^{2} \equiv 0 \bmod p$. 

Similarly, if $q$ is an odd prime then $V_{q} \equiv qh^{q-1} \bmod dg^{2}$. 
Since $p|(dg^{2})$, we have  
\begin{equation}
\label{eq:bqcong}
V_{q} \equiv qh^{q-1} \bmod p. 
\end{equation}

Notice that $V_{q} \not\equiv 0 \bmod p$, for if this were true 
then $p|h$, but $p|(dg^{2})$, so $p|(h^{2}-dg^{2})$, which is 
not possible since $h^{2}-dg^{2}=-1$.  

We are now ready to show that our assumption that $V_{p}=py_{1}^{2}$ 
is impossible, from which our lemma will follow. 

We can apply Lemma~\ref{lem:2.11}(iii), since $\alpha+\beta=2h$. 
If $V_{p}=py_{1}^{2}$ then $pV_{p}$ would be a square and, for 
any odd prime $q \neq p$,  
\begin{displaymath}
1 = \left( \frac{pV_{p}}{V_{q}} \right) 
  = \left( \frac{p}{V_{q}} \right) \left( \frac{V_{p}}{V_{q}} \right)   
  = \left( \frac{p}{V_{q}} \right) 
  = \left( \frac{V_{q}}{p} \right),   
\end{displaymath}  
the third equality holding by Lemma~\ref{lem:2.11}(iii) 
and the last one by the law of quadratic reciprocity, 
since $V_{q} \equiv 1 \bmod 4$, as shown in the proof 
of Lemma~\ref{lem:2.11}(iii). 

Applying (\ref{eq:bqcong}), we then have 
\begin{displaymath}
1 = \left( \frac{qh^{q-1}}{p} \right) 
  = \left( \frac{q}{p} \right).  
\end{displaymath}  

But clearly this cannot hold for all odd primes $q \neq p$ and so 
we obtain our desired contradiction. Therefore $V_{p} \neq py_{1}$ 
and $t$, the least positive integer such that (\ref{eq:chena}), 
must be one. The lemma follows. 
\end{proof}

\vspace{3.0mm}

\begin{center}
3. {\sc Proof Of Theorem~\ref{thm:E}} 
\end{center}

\setcounter{equation}{0}
\setcounter{lemma}{0}
\addtocounter{section}{1}

\vspace{3.0mm}

We shall use the hypergeometric method to solve 
(\ref{eq:3.1}) for $t \geq 128$. One could refine 
the arguments in the previous section to solve 
some of the remaining Thue equations but not all 
of them: some further method is needed. This is 
the same situation which arises in \cite{Thom} 
and similar papers, including \cite{LP}. Since 
Lettl and Peth\H{o} have considered the remaining 
equations in their paper by the same method we would 
wish to use here, it seems repetitive to produce an 
argument here showing that Theorem~\ref{thm:E} holds 
for $t \leq 127$. Therefore, we shall refer the reader 
to \cite{LP} for $1 \leq t \leq 127$ and only prove 
this theorem for $t \geq 128$ in what follows. 

For any positive integer $t$, let us define  
\begin{eqnarray*}
\epsilon = \frac{t + \sqrt{t^{2}+16}}{4}, 
& & \rho = \sqrt{1+\epsilon^{2}}, \\ 
\beta^{(0)} = \epsilon - \rho, \hspace{5 mm}  
\beta^{(1)} = \frac{\rho-1}{\epsilon}, 
& & \beta^{(2)} = \rho+\epsilon 
    \hspace{3.0mm} \mbox{ and } \hspace{3.0mm}
    \beta^{(3)} = -\frac{\rho+1}{\epsilon}. 
\end{eqnarray*}

One can show that the $\beta^{(j)}$'s are the four roots of $P_{t}(X,1)$.  

We start by using the results from Section~2 to construct an infinite 
sequence of rational approximations to the $\beta^{(j)}$'s which 
are sufficiently good to allow us to obtain an effective improvement 
on Liouville's theorem for these $\beta^{(j)}$'s. We will use this 
result to solve (\ref{eq:3.1}).

\vspace{3.0mm}

\begin{center}
3.1. {\sc Rational Approximations To The $\beta^{(j)}$'s} 
\end{center}

\vspace{3.0mm}

We shall construct these rational approximations from Lemma~\ref{lem:2.3}. 
Since $\beta^{(0)} \beta^{(2)}=\beta^{(1)} \beta^{(3)} =-1$, if 
$p_{0}-\beta^{(0)}q_{0}$ and $p_{1}-\beta^{(1)}q_{1}$ are close to 
zero then so are $p_{0}\beta^{(2)}+q_{0}$ and $p_{1}\beta^{(3)}+q_{1}$. 
Therefore, we will initially  only consider approximations to 
$\beta^{(0)}$ and $\beta^{(1)}$. 

In the remainder of this section, we shall assume that $t$ is a fixed 
integer greater than 127. 

We now determine the quantities defined in the Lemma~\ref{lem:2.3}. 
Put $P(X) = X^{4}-tX^{3}-6X^{2}+tX+1$ and $U(X)=X^{2}+1$. It is easy 
to check that $P(X)$ and $U(X)$ satisfy the differential equation 
(\ref{eq:2.8}) and that the discriminant of $U(X)$ is non-zero. 
Therefore, Lemma~\ref{lem:2.3} is applicable. 

Other important quantities are 
\begin{eqnarray*}
h = -15, \hspace{5 mm} \lambda = -1, & & 
Y_{1}(X) = 2tX^{4}+32X^{3}-12tX^{2}-32X+2t, \\ 
A_{0}(X) = -10, & & 
A_{1}(X) = \frac{10X^{5}}{3}+\frac{100X^{3}}{3}-\frac{50tX^{2}}{3}
	   -50X+\frac{10t}{3}, \\ 
B_{0}(X) = -10X, & & 
B_{1}(X) = \frac{10tX^{5}}{3}+50X^{4}-\frac{50tX^{3}}{3}
	   -\frac{100X^{2}}{3}-\frac{10}{3}, \\ 
z(X) = \frac{1}{16} (-Y_{1}(X)i + 8P(X)), & &
u(X) = \frac{1}{16} (-Y_{1}(X)i - 8P(X)), \\ 
a(X) =  5iX-5, &            & b(X) = 5iX+5, \\  
c(X) = -5X-5i  & \mbox{and} & d(X) = 5X-5i. 
\end{eqnarray*}

We could use Lemma~\ref{lem:2.3}(i) with any rational $x$ 
to construct a sequence of rational approximations, but 
Lemma~\ref{lem:2.5} suggests that choosing $x$ so that 
$\beta^{(j)} (a(x)w(x)^{1/n}-b(x))-(c(x)w(x)^{1/n}-d(x))=0$ 
would give us a sequence of particularly good approximations 
to $\beta^{(j)}$. 

Fortunately, we can find such values. For $j=0$, we can choose 
$x=0$ and for $j=1$, we can choose $x=1$. To see this notice 
that 
\begin{eqnarray*}
P(0) = 1, \hspace{5 mm} P(1) = -4, \hspace{5 mm} 
Y_{1}(0) = 2t \mbox{\hspace{3.0mm} and \hspace{3.0mm}} Y_{1}(1) = -8t. 
\end{eqnarray*}

Hence 
\begin{eqnarray*}
u(0) = \frac{-it-4}{8}, \hspace{3.0mm} 
z(0) = \frac{-it+4}{8}, \hspace{3.0mm} 
u(1) = \frac{it+4}{2} \hspace{3.0mm} \mbox{and} \hspace{3.0mm}
z(1) = \frac{it-4}{2}. 
\end{eqnarray*}

So we have 
\begin{displaymath}
w(0) = w(1) = \frac{z(0)}{u(0)} = \frac{z(1)}{u(1)} 
= \frac{t^{2}-16}{t^{2}+16} + i \frac{8t}{t^{2}+16}.   
\end{displaymath}

In what follows, we let $w$ denote this number. For 
$t > 4$, we can write $w = e^{i \varphi}$, where 
$0 < \varphi < \pi/2$. 

Notice that 
\begin{displaymath}
\sqrt{w} = e^{i \varphi/2} = \frac{t+4i}{\sqrt{t^{2}+16}} 
\mbox{ and } 
\sqrt[4]{w} = \frac{\epsilon}{\rho} + i \frac{1}{\rho}, 
\end{displaymath}
where $\epsilon$ and $\rho$ are as defined at the beginning 
of Section~3. 

Finally, we have 
\begin{eqnarray*}
a(0) = -5,   \hspace{9.0mm} b(0) = 5, \hspace{7.0mm} & & c(0) = -5i, 
\hspace{12.0mm} d(0) = -5i \nonumber \\
a(1) = 5i-5, \hspace{5 mm} b(1) = 5i+5,              & & c(1) = -5-5i   
\mbox{ and } d(1) = 5-5i.        
\end{eqnarray*}

So now it is a routine matter to verify that, for $j=0$ and 1,   
\begin{displaymath}
\beta^{(j)} \left( a(j) w^{1/4} - b(j) \right) 
- \left( c(j) w^{1/4} - d(j) \right) = 0,  
\end{displaymath}
so the first term in the expression for $i^{r}C_{r}(0)$ and 
$i^{r}C_{r}(1)$ in Lemma~\ref{lem:2.5} disappears. 

\vspace{3.0mm}

\noindent
{\bf Notation.} We shall simplify our notation further here to reflect 
the fact that we are only considering $n=4$. We shall use $R_{r}$ and 
$X_{r}$ instead of $R_{4,r}$ and $X_{4,r}$.  

\vspace{3.0mm}

We now construct our sequences of rational approximations to 
$\beta^{(0)}$ and $\beta^{(1)}$. 

From Lemmas~\ref{lem:2.3} and \ref{lem:2.5}, we have, for $j=0,1$, 
\begin{eqnarray}
\label{eq:3.4}
i^{r} A_{r}(j) & = & a(j) X_{r}^{*}(z(j),u(j)) 
		       - b(j) X_{r}^{*}(u(j),z(j)),  \nonumber \\ 
i^{r} B_{r}(j) & = & c(j) X_{r}^{*}(z(j),u(j)) 
		       - d(j) X_{r}^{*}(u(j),z(j)) \mbox{ and } \\ 
i^{r} C_{r}(j) & = & -\left( \beta^{(j)} a(j) - c(j) \right) 
		       u(j)^{r} R_{r}(w). \nonumber
\end{eqnarray}

These quantities will form the basis for our approximations. 
We first eliminate some common factors. 
Putting $D_{r}(0)=8^{r}/5$ and $D_{r}(1)=(-2)^{r}/5$, 
we can now define 
\begin{displaymath}
P_{r}^{(j) \prime} = D_{r}(j)B_{r}(j) 
\hspace{3.0mm} \mbox{ and } \hspace{3.0mm}
Q_{r}^{(j) \prime} = D_{r}(j)A_{r}(j), 
\end{displaymath}
for $j=0$ and 1. 

Furthermore, defining $z' = -it+4$ and $u' = -it-4$,  
we have  
\begin{displaymath}
i^{r} P_{r}^{(j) \prime} = \frac{c(j)X_{r}^{*}(z',u') - d(j)X_{r}^{*}(u',z')}
				{5}  
\end{displaymath}
and 
\begin{displaymath}
i^{r} Q_{r}^{(j) \prime} = \frac{a(j)X_{r}^{*}(z',u') - b(j)X_{r}^{*}(u',z')}
				{5}. 
\end{displaymath}

Siegel has shown \cite[Hilfssatz 4 and the proof of Hilfssatz 3]{Siegel} that  
we can write 
\begin{eqnarray*}
X_{r}^{*}(u',z') & = & \frac{4^{r} (r!)}{3 \cdot 7 \cdots (4r-1)} 
		       {2r \choose r} 
		       (z')^{r} {} _{2}F_{1} \left( -r,-r-1/4,-2r,1-w^{-1} \right)
		       \mbox{ and } \\ 
X_{r}^{*}(z',u') & = & \frac{4^{r} (r!)}{3 \cdot 7 \cdots (4r-1)} 
		       {2r \choose r} 
		       (u')^{r} {} _{2}F_{1}(-r,-r-1/4,-2r,1-w).  
\end{eqnarray*}

Notice that $1-w=-8/u', 1-w^{-1}=8/z'$ and $\mu_{4}=2$ so that, by 
Lemma~\ref{lem:2.6}, $MX_{r}^{*}(u',z')$ and $MX_{r}^{*}(z',u')$ 
are algebraic integers where 
\begin{displaymath}
M = \frac{3 \cdot 7 \cdots (4r-1)}{4^{r} r!}. 
\end{displaymath}

This implies that $MP_{r}^{(j) \prime}$ and $MQ_{r}^{(j) \prime}$ 
are also algebraic integers. In addition, we have the relation 
$z'=-\overline{u'}$ so that $X_{r}^{*}(u',z') 
= \pm \overline{X_{r}^{*}(z',u')}$. Thus 
$P_{r}^{(j)}=MP_{r}^{(j) \prime}/2$ and 
$Q_{r}^{(j)}=MQ_{r}^{(j) \prime}/2$ are algebraic integers too. 

Moreover, since $P(X), U(X), Y_{1}(X) \in \bbZ[X]$ and 
$h, \lambda \in \bbZ$, $A_{r}(X)$ and $B_{r}(X)$ have 
rational coefficients, by their definition in Lemma~\ref{lem:2.3}, 
and so $P_{r}^{(j) \prime}$ and $Q_{r}^{(j) \prime}$ belong 
to $\bbQ$. Therefore, $P_{r}^{(j)}$ and $Q_{r}^{(j)}$ are, 
in fact, rational integers.  

These are the numbers we shall use for our rational approximations. 
We have 
\begin{displaymath}
Q_{r}^{(j)} \beta^{(j)} - P_{r}^{(j)} = S_{r}^{(j)}, 
\end{displaymath}
for $j=0$ and 1 where $S_{r}^{(j)}=D_{r}(j)MC_{r}(j)/2$. 

We now want to show that these are ``good'' approximations; 
we do this by estimating $|P_{r}^{(j)}|, |Q_{r}^{(j)}|$ and 
$|S_{r}^{(j)}|$ from above. 

We can write $M$ as 
\begin{displaymath}
\frac{\Gamma(r+3/4)}{\Gamma(3/4) \, r!}. 
\end{displaymath}

Note that $|a(0)|=|b(0)|=|c(0)|=|d(0)|=5,
|a(1)|=|b(1)|=|c(1)|=|d(1)|=5\sqrt{2}, |u'|=|z'|$ and 
$| u'(1+\sqrt{w})^{2}| = 8\epsilon$. By Lemma~\ref{lem:2.8} 
and the triangle inequality, we see that 
\begin{eqnarray*} 
|P_{r}^{(0)}| = \frac{M |P_{r}^{(0) \prime}|}{2} 
& \leq & \frac{\Gamma(r+3/4)}{\Gamma(3/4) \, r!} 
	 \cdot \frac{4}{|1+\sqrt{w}|^{2}} 
	 \cdot \frac{r! \, \Gamma(3/4)}{\Gamma(r+3/4)} 
	 { \left| u'(1+\sqrt{w})^{2} \right| }^{r} \nonumber \\ 
& = &    \frac{4}{|1+\sqrt{w}|^{2}} 
	 { \left| u'(1+\sqrt{w})^{2} \right| }^{r} 
=        \frac{4}{|1+\sqrt{w}|^{2}} (8\epsilon)^{r}.   
\end{eqnarray*} 

Since $|1+\sqrt{w}|^{2} > 3.999$ for $t \geq 128$, we have 
\begin{equation}
\label{eq:3.12}
|P_{r}^{(0)}| < 1.0005 (8\epsilon)^{r}. 
\end{equation}

Similarly, we have 
\begin{eqnarray*}
|P_{r}^{(1)}| & < & 1.415 (8\epsilon)^{r}, \\ 
|Q_{r}^{(0)}| & < & 1.0005 (8\epsilon)^{r} \mbox{ and } \\  
|Q_{r}^{(1)}| & < & 1.415 (8\epsilon)^{r},  
\end{eqnarray*}
since $1.0005 \sqrt{2} < 1.415$.

By Lemma~\ref{lem:2.1} and (\ref{eq:3.4}), we have, for $j=0,1$,  
\begin{displaymath}
|S_{r}^{(j)}| \leq \left| D_{r}(j) \right| 
		   \frac{\Gamma(r+3/4)}{2\Gamma(3/4) \, r!} 
		   \left| \beta^{(j)} a(j)-c(j) \right| |u(j)|^{r}
		   \frac{\Gamma(r+5/4)}{\Gamma(1/4) \, r!} 
		   \varphi { \left| 1-\sqrt{w} \right| }^{2r}.  
\end{displaymath}

Since 
\begin{eqnarray*}
\left| D_{r}(j) \right| \left| \beta^{(j)} a(j)-c(j) \right| |u(j)|^{r} 
& = & \frac{\left| \beta^{(j)} a(j)-c(j) \right| |u'|^{r}}{5}, \\  
\frac{2 \varphi}{\pi} \leq \sin \varphi & = & \frac{8t}{16+t^{2}}, \\ 
\frac{\Gamma(r+3/4)}{2\Gamma(3/4) \, r!} 
\frac{\Gamma(r+5/4)}{\Gamma(1/4) \, r!} 
& \leq & \frac{1}{8} \mbox{ for $r \geq 0$ and } \\
\left| u'(1-\sqrt{w})^{2} \right| & = & \frac{8}{\epsilon}, 
\end{eqnarray*}
we have 
\begin{displaymath}
|S_{r}^{(j)}| \leq \frac{\pi}{2} \left| \beta^{(j)}a(j) - c(j) \right|  
		   \frac{t}{16+t^{2}} 
		   { \left( \frac{8}{\epsilon} \right) }^{r}. 
\end{displaymath}

Since $|\beta^{(j)}a(j) - c(j)| \leq 10$, we obtain 
\begin{equation}
\label{eq:3.13}
|S_{r}^{(j)}| \leq \frac{\pi t}{16+t^{2}} 
	     { \left( \frac{8}{\epsilon} \right) }^{r}, 
\end{equation}
for $j=0$ and 1.

Noting that $\beta^{(0)}\beta^{(2)}=\beta^{(1)}\beta^{(3)}=-1$, 
we have 
\begin{eqnarray}
\label{eq:3.14}
-Q_{r}^{(0)}-\beta^{(2)}P_{r}^{(0)} & = & \beta^{(2)} S_{r}^{(0)}
\hspace{3.0mm} \mbox{ and } \\ 
\label{eq:3.15}
-Q_{r}^{(1)}-\beta^{(3)}P_{r}^{(1)} & = & \beta^{(3)} S_{r}^{(1)}.
\end{eqnarray}

Now we apply Lemma~\ref{lem:2.9} to prove the following theorem.

\begin{theorem}
%\label{thm:D}
Let the $\beta^{(j)}$'s and $\epsilon$ be as above and 
suppose that $t \geq 128$ is a rational integer. Define 
\begin{displaymath}
\kappa = \frac{\log 8\epsilon}{\log \epsilon/8}.
\end{displaymath}

For $j=0,1,2,3$ and any rational integers $p$ and $q$, we have  
\begin{equation}
\label{eq:3.16}
\left| p - \beta^{(j)} q \right| > \frac{1}{c_{j}|q|^{\kappa}} 
\end{equation}
for $|q| > 0.16t$, where  
\begin{displaymath}
c_{0} = c_{1} = c_{3} = 11.33t \cdot 0.4^{\kappa}, \hspace{3.0mm} 
\mbox{ and } \hspace{3.0mm} c_{2} = 8.01t (0.4t)^{\kappa}.  
\end{displaymath}
\end{theorem}

\begin{proof}
In each case we will apply Lemmas~\ref{lem:2.7} and \ref{lem:2.9}. 

We have 
\begin{displaymath}
P_{r}^{(j)}Q_{r+1}^{(j)}-P_{r+1}^{(j)}Q_{r}^{(j)}
= { \left( \frac{M}{2} \right) }^{2} D_{r}(j) D_{r+1}(j)
  \left( A_{r+1}(j)B_{r}(j)-A_{r}(j)B_{r+1}(j) \right). 
\end{displaymath}

Applying Lemma~\ref{lem:2.7} with $a=d=1, b=c=0$ and $X=j$, we 
see that $P_{r}^{(j)}Q_{r+1}^{(j)} \neq P_{r+1}^{(j)}Q_{r}^{(j)}$.

Letting $p_{r}=P_{r}^{(0)}$ and $q_{r}=Q_{r}^{(0)}$, from (\ref{eq:3.12}) 
and (\ref{eq:3.13}), we can take $k_{0}=1.0005, l_{0} = \pi t/(16+t^{2}), 
E=\epsilon/8$ and $Q=8\epsilon$. Therefore, we see that $2l_{0}E < 0.4$ 
and $2k_{0}Q < 16.008 \epsilon < 8.01t$ for $t \geq 128$. Hence we can 
use $c_{0}$ for the quantity $c$ in Lemma~\ref{lem:2.9}. For $t \geq 128$, 
$1/(2l_{0})=(16+t^{2})/(2\pi t) < 0.16t$. 
This gives us the lower bound $N_{0}$ for $|q|$. 

For $j=1$, the use of Lemma~\ref{lem:2.9} is identical to 
its use for $j=0$ except that here $k_{0}=1.415$.

By (\ref{eq:3.14}), for $j=2$, we simply switch $p_{r}$ and $q_{r}$ 
from $j=0$ and let $l_{0}$ is $\beta^{(2)}$ times the value of 
$l_{0}$ used for $j=0$. One can show that $2l_{0}E < 0.4t$ 
and $1/(2l_{0}) < 0.16$ for $t \geq 128$ which implies that we 
can use $c_{2}$ for $c$ and $N_{2}$ as a lower bound for $|q|$ 
in Lemma~\ref{lem:2.9}. 

Finally, by (\ref{eq:3.15}), the theorem for $j=3$ follows from 
the result for $j=1$ and the inequality $\epsilon/(\rho+1) < 1$, 
upon switching $p_{r}$ and $q_{r}$.  
\end{proof}

\begin{center}
3.2. {\sc Solving The Thue Equation (\ref{eq:3.1}) For $t \geq 128$} 
\end{center}

\vspace{3.0mm}

As we just saw, the hypergeometric method will give us an  
irrationality measure for the $\beta^{(j)}$'s once $q$ is 
sufficiently large. So let us first deal with $|q|$ small. 

To get started, we want to show that if $(x,y)$ is a solution 
of (\ref{eq:3.1}) with $|y| > 1$ then $x/y$ is a convergent 
in the continued-fraction expansion to one of the $\beta^{(j)}$'s. 
For this we need to know that, for $t \geq 5$,  
\begin{eqnarray*}
- \frac{1}{t}      < & \beta^{(0)} & < - \frac{1}{t+1}, \\  
1 - \frac{2}{t+1}  < & \beta^{(1)} & < 1 - \frac{2}{t+2}, \\ 
t + \frac{5}{t+1}  < & \beta^{(2)} & < t + \frac{5}{t} \mbox{ and } \\
-1 - \frac{2}{t-1} < & \beta^{(3)} & < -1 - \frac{2}{t}.  
\end{eqnarray*}

From these and some calculation for small $t$, one can verify that,  
for $t \geq 1$, 
\begin{eqnarray*}
t           < & f'(\beta^{(0)}) & < t + \frac{9}{t}, \\ 
-2t - 8     < & f'(\beta^{(1)}) & < -2t+4, \\ 
t^{3} + 19t < & f'(\beta^{(2)}) & < t^{3}+38t \mbox{ and } \\  
-2t - 20    < & f'(\beta^{(3)}) & < -2t - 4. 
\end{eqnarray*}

Therefore, $C_{1} = \min_{j} |f'(\beta^{(j)})| > 32$ for $t \geq 32$. 

As a particular case of Lemma~1.1 of \cite{TW}, we know that  
if $|y| > 1 \geq \lceil \sqrt{32/C_{1}} \rceil$ then $x/y$ is 
a convergent to one of the $\beta^{(j)}$'s. 

Let us first consider $y=0$ and $|y|=1$. 

If $y=0$ then $P_{t}(x,0)=x^{4}=\pm 1$ implies $x= \pm 1$. 

If $y=1$ then $P_{t}(x,1)=x^{4}-tx^{3}-6x^{2}+tx+1=1$ implies 
that $x(x^{3}-tx^{2}-6x+t)=0$. For $1 \leq t \leq 5$, the only 
rational roots of this polynomial in $x$ is $x=0$. For $t \geq 5$, 
we have either $x=0, -2 < x < -1, 0 < x < 1$ or $t < x < t+1$. 

If $P_{t}(x,1)=x^{4}-tx^{3}-6x^{2}+tx+1=-1$ then 
$x^{4}-tx^{3}-6x^{2}+tx+2=0$ which has the rational root $x=-2$ 
for $t=1$ and no others for $1 \leq t \leq 5$. For $t \geq 5$, 
we have either $-2 < x < -1, -1 < x < 0, 0 < x < 1$ or $t < x < t+1$. 

If $y=-1$ then $P_{t}(x,-1)=1$ implies that $x(x^{3}+tx^{2}-6x-t)=0$. 
For $1 \leq t \leq 5$, the only rational roots of this polynomial 
in $x$ is $x=0$. For $t \geq 5$, we have either $-(t+1) < x < -t, 
-1 < x < 0$ or $1 < x < 2$. 

If $P_{t}(x,-1)=-1$ then there is the rational root $x=-2$ 
for $t=1$ and no others for $1 \leq t \leq 5$. For $t \geq 5$, 
have either $-(t+1) < x < -t, -1 < x < 0, 0 < x < 1$ or $1 < x < 2$. 

So we see that for $y=0$ or $|y|=1$, there are no solutions of 
(\ref{eq:3.1}) except those mentioned in Theorem~\ref{thm:E}. 

Using the bounds above for $\beta^{(0)}$, we find that the 
continued-fraction expansion for $\beta^{(0)}$ is $[-1,1,t-1,\ldots]$ 
when $t \geq 5$. Therefore the convergents of $\beta^{(0)}$ 
are $-1/1,0/1,-1/t,\ldots$ for $t \geq 5$. In the same manner, 
we find that, for $t \geq 5$, 
$\beta^{(1)} = [ 0, 1, \lfloor (t-1)/2 \rfloor, \ldots ], 
\beta^{(2)} = [ t, \lfloor t/5 \rfloor, \ldots ]$ and 
$\beta^{(3)} = [ -2, 1, \lfloor (t-3)/2 \rfloor, \ldots ]$ and 
their convergents are, respectively, 
$0/1,1/1,\lfloor (t-1)/2 \rfloor / (\lfloor (t-1)/2 \rfloor + 1), \ldots$ 
and $t/1, (\lfloor t/5 \rfloor t + 1) / \lfloor t/5 \rfloor, \ldots$ and  
$-2/1,-1/1,-(\lfloor (t-3)/2 \rfloor + 2) / (\lfloor (t-3)/2 \rfloor + 1), 
\ldots$. 

Notice that the smallest of these denominators which is not 1 is 
at least $(t-4)/5$, which is greater than $0.193t$ for $t \geq 128$. 
Hence (\ref{eq:3.1}) has no solutions with $1 < |y| < 0.193t$ for 
$t \geq 128$. 

Now we can turn to the `large' solutions. 

Let $\delta^{(j)}=|x-y\beta^{(j)}|$ for $0 \leq j \leq 3$, 
where the $\beta^{(j)}$'s are as above. If $(x,y)$ 
is a solution of (\ref{eq:3.1}) then 
\begin{equation}
\label{eq:3.19}
\delta^{(0)} \delta^{(1)} \delta^{(2)} \delta^{(3)} = 1. 
\end{equation}

If $y \neq 0$ then the $\delta^{(j)}$'s distinct and there 
exists a smallest one which is less than 1. 

If $\delta^{(0)}$ is the smallest of the $\delta^{(j)}$'s then 
\begin{equation}
\label{eq:3.19a}
\delta^{(j)} = |x- \beta^{(0)}y+(\beta^{(0)}-\beta^{(j)})y| 
> |\beta^{(0)}-\beta^{(j)}| |y| - 1 \mbox{ for $j=1,2$ and 3.}
\end{equation}

For $t \geq 128$, it is easy to show from (\ref{eq:3.19a}) 
and the estimates at the beginning of Section~3.2 that we have 
\begin{displaymath}        
\delta^{(1)} > 0.99 |y| - 1, \hspace{5 mm} 
\delta^{(2)} > t|y|-1 \mbox{ and }
\delta^{(3)} > |y| - 1. 
\end{displaymath}

From (\ref{eq:3.19}) and the above inequalities, we see that 
\begin{displaymath}
\delta^{(0)} < \frac{1}{|y|^{3}(1-1/|y|)(0.99-1/|y|)(t-1/|y|)} 
< \frac{1}{0.911t|y|^{3}},  
\end{displaymath}
for $|y| \geq \lfloor t/5 \rfloor \geq 25$.

Combining this inequality with (\ref{eq:3.16}), we see that 
\begin{displaymath}
|y|^{3-\kappa} 
< 12.5 (0.4)^{\kappa}. 
\end{displaymath}

For $t \geq 128$, we have $1 < \kappa < 3$ and $|y| \geq 25$, 
so $25^{3-\kappa} \leq |y|^{3-\kappa}$. But for $\kappa$ in 
this range, it is easy to see that $12.5 (0.4)^{\kappa} 
< 25^{3-\kappa}$, so the required inequality is not true. 
Therefore, there are no large solutions of (\ref{eq:3.1}) for 
which $\delta^{(0)}$ is the smallest of the $\delta^{(j)}$'s 
when $t \geq 128$. 

With a similar method, we can show that the same is true 
for the other there are no large $\delta^{(j)}$'s.  

If $\delta^{(1)}$ is the smallest then 
\begin{displaymath}
\delta^{(1)} < \frac{1}{1.84t|y|^{3}}, 
\end{displaymath}
for $t \geq 128$ and $y \geq \lfloor t/5 \rfloor \geq 25$. 

We derive a contradiction in the same way as in the case of $\delta^{(0)}$. 

If $\delta^{(2)}$ is the smallest of the $\delta^{(j)}$'s then 
\begin{displaymath}
\delta^{(2)} < \frac{1}{0.999t^{3}|y|^{3}}, 
\end{displaymath}
for $t \geq 128$ and $y \geq \lfloor t/5 \rfloor \geq 25$. 

Combining this inequality with (\ref{eq:3.16}), we get 
\begin{displaymath}
|y|^{3-\kappa} 
< 8.02t^{-2}(0.4t)^{\kappa}. 
\end{displaymath}

Since $0.19t < \lfloor t/5 \rfloor \geq |y|$ for $t \geq 128$, 
we want $t^{3} < 8.02 (0.19)^{\kappa-3}(0.4)^{\kappa} < 89$ 
since $1 < \kappa$. But this implies that $t < 5$ --- a 
contradiction. 

Finally, if $\delta^{(3)}$ is the smallest then 
\begin{displaymath}
\delta^{(3)} < \frac{1}{1.89t|y|^{3}} 
\end{displaymath}
and we obtain a contradiction as with $\delta^{(0)}$. 

Thus there are no large solutions to (\ref{eq:3.1}) either 
and Theorem~\ref{thm:E} follows for $t \geq 128$. 

\vspace{3.0mm}

\begin{center}
4. {\sc Proof Of Theorem~\ref{thm:C}} 
\end{center}

\setcounter{equation}{0}
\setcounter{lemma}{0}
\addtocounter{section}{1}

\vspace{3.0mm}

If there are no solutions to (\ref{eq:1.1}) then we are done, 
so let us suppose otherwise. 

By Lemma~\ref{lem:2.12}, $(x_{0},y_{0})$ is a solution of 
(\ref{eq:1.1}) where $(x_{0},y_{0}^{2})$ is the fundamental 
solution of the Pell equation $X^{2}+1=dY^{2}$. 

Putting $\epsilon = x_{0}+y_{0}^{2}\sqrt{d}$, since 
$x_{0}^{2}+1=dy_{0}^{4}$, we have $y_{0}^{2}\sqrt{d} 
= \sqrt{1+x_{0}^{2}}$ and hence 
\begin{displaymath}
\epsilon = x_{0} + \sqrt{x_{0}^{2}+1}. 
\end{displaymath}

If (\ref{eq:1.1}) has another solution $(x,y)$, then we have 
\begin{displaymath}
y^{2} = \frac{\epsilon^{2t+1} - {\overline{\epsilon}}^{2t+1}}
	     {2 \sqrt{d}},  
\end{displaymath}             
where $t \in \bbZ$ and $\overline{\epsilon}=x_{0}-\sqrt{x_{0}^{2}+1}$. 

Hence we have 
\begin{displaymath}
{ \left( \frac{y}{y_{0}} \right) }^{2} 
= \frac{\epsilon^{2t+1} - {\overline{\epsilon}}^{2t+1}}{2 y_{0}^{2}\sqrt{d}} 
= \frac{\epsilon^{2t+1} - {\overline{\epsilon}}^{2t+1}}{2 \sqrt{1+x_{0}^{2}}}.  
\end{displaymath}

Note that 
\begin{displaymath}
\epsilon^{2t+1} = U_{2t+1}+V_{2t+1}\sqrt{1+x_{0}^{2}} 
\hspace{5 mm} \mbox{ and } \hspace{5 mm} 
\overline{\epsilon}^{2t+1} = U_{2t+1}-V_{2t+1}\sqrt{1+x_{0}^{2}}, 
\end{displaymath}
where $U_{2t+1}, V_{2t+1} \in \bbZ$. Therefore, 
$(y/y_ {0})^{2} = V_{2t+1} \in \bbZ$, so $y/y_{0}$ is an algebraic 
integer and hence a rational integer. 

The above analysis shows that if $(x,y)$ is a solution of (\ref{eq:1.1}) 
then $(x,y/y_{0})$ is a solution of 
\begin{equation}
\label{eq:4.1}
X^{2}+1=(1+x_{0}^{2})Y^{4}. 
\end{equation}

Conversely, given $(x,y) \in \bbZ^{2}$ which satisfy 
$x^{2}+1=(1+x_{0}^{2})y^{4}$ for some $x_{0} \in \bbZ$ 
and a $d \in \bbZ$ such that $1+x_{0}^{2}=dy_{0}^{4}$ 
for some $y_{0} \in \bbZ$, then $(x,yy_{0})$ is a solution 
of (\ref{eq:1.1}). So we only need to consider the case 
$d=1+x_{0}^{2}$ and in the remainder of this section, 
we only look at the equation (\ref{eq:4.1}). 

The equation (\ref{eq:4.1}) has an obvious solution, 
namely $(x_{0},1)$, which corresponds to the fundamental 
solution of $X^{2}+1=(1+x_{0}^{2})Y^{2}$. 

As above, we write 
\begin{displaymath}
\epsilon^{m} = U_{m}+V_{m}\sqrt{1+x_{0}^{2}},  
\end{displaymath}
for non-negative integers $m$. Notice that $U_{1}=x_{0}$ and $V_{1}=1$. 

As before, if (\ref{eq:4.1}) has another solution $(x,y)$, then 
\begin{equation}
\label{eq:4.3}
y^{2} = \frac{\epsilon^{2t+1}-\overline{\epsilon}^{2t+1}}
	     {2\sqrt{1+x_{0}^{2}}} 
      = V_{2t+1}, 
\end{equation} 
for some non-negative integer $t$. 

It is easy to see that 
\begin{eqnarray}
\label{eq:4.4}
V_{m+1}  & = & x_{0}V_{m}+U_{m} \hspace{3.0mm} \mbox{ for $m \geq 1$ and } \\
\label{eq:4.5}
V_{2m+1} & = & V_{1} V_{2m+1} = V_{m}^{2} + V_{m+1}^{2}.
\end{eqnarray}

Combining (\ref{eq:4.3}) with (\ref{eq:4.5}), we get 
\begin{equation}
\label{eq:4.6}
y^{2}=V_{t}^{2}+V_{t+1}^{2}. 
\end{equation}

We distinguish two cases. 

If $2 | V_{t}$ then, from (\ref{eq:4.6}) and the 
theory of Pythagorean numbers, we have 
\begin{equation}
\label{eq:4.7}
V_{t} = 2ab, \hspace{5.0mm} V_{t+1} = a^{2}-b^{2} 
\hspace{3.0mm} \mbox{ and } \hspace{3.0mm} 
y = \pm \left( a^{2}+b^{2} \right),  
\end{equation}
where $a,b \in \bbZ$ and $(a,b)=1$, this latter condition 
holding because of Lemma~\ref{lem:2.11}(ii). 

From (\ref{eq:4.4}) and (\ref{eq:4.7}), we get  
\begin{equation}
\label{eq:4.8}
U_{t} = V_{t+1}-x_{0}V_{t} = a^{2}-b^{2}-2x_{0}ab. 
\end{equation}

By the definition of $U_{t}$ and $V_{t}$, 
\begin{displaymath}
\epsilon^{t} \overline{\epsilon}^{t}
= U_{t}^{2}-(1+x_{0}^{2})V_{t}^{2} = \pm 1.  
\end{displaymath}

So, by (\ref{eq:4.7}) and (\ref{eq:4.8}), we find that 
\begin{displaymath}
(a^{2}-b^{2}-2x_{0}ab)^{2} - 4(x_{0}^{2}+1)a^{2}b^{2} = \pm 1.  
\end{displaymath}

Hence we have 
\begin{displaymath}
a^{4}-4x_{0}a^{3}b-6a^{2}b^{2}+4x_{0}ab^{3}+b^{4} = \pm 1.  
\end{displaymath}

Applying Theorem~\ref{thm:E} with $t=4x_{0}$, we see that 
$y=\pm(a^{2}+b^{2})=\pm 1$ and the theorem follows. 

If 2 does not divide $V_{t}$ then considering 
(\ref{eq:4.6}) modulo 4, we see that $2 | V_{t+1}$. 
So 
\begin{displaymath}
V_{t} = a^{2}-b^{2}, \hspace{5.0mm} V_{t+1} = 2ab 
\hspace{3.0mm} \mbox{ and } \hspace{3.0mm} 
y = \pm \left( a^{2}+b^{2} \right),  
\end{displaymath}
where $a,b \in \bbZ$ and $(a,b)=1$. 

Proceeding as in the previous case, we obtain the equation
\begin{displaymath}
(-a)^{4}-4x_{0}(-a)^{3}b-6(-a)^{2}b^{2}+4x_{0}(-a)b^{3}+b^{4} = \pm 1.  
\end{displaymath}

Again, the theorem follows from Theorem~\ref{thm:E}. 

\vspace{1.0mm}

\begin{center}
{\sc Acknowledgements}
\end{center}

\vspace{1.0mm}

Both authors would like to express their deepest thanks to Professor 
M. Waldschmidt for, amongst many things, making their individual work 
known to the other. In addition, the first author wishes to express 
his gratitude to Professors M. Mignotte, Qi Mingyou and Zeng Xianwu 
for their kind help as well as to his wife Liu Ying for her assistance. 

Finally, the authors would like to thank the referee for his/her 
careful reading of our manuscript and the suggestions made, especially 
regarding Lemma~\ref{lem:2.9}. 

2018 addendum: thanks as well to Bernadette Faye and Eva Goedhart for their
very attentive reading, finding that our proof violated an assumption made
in Lemma~\ref{lem:2.5}.

\end{document}